%% file: rigid.tex
\documentclass{article}
\usepackage{amsmath,amssymb,theorem,xypic}
\xyoption{curve}

\theoremstyle{change}
\theorembodyfont{\itshape}
\newtheorem{thm}{Theorem}
\newtheorem{prop}{Proposition}
\newtheorem{lemma}{Lemma}

\numberwithin{equation}{section}

\newcommand{\demobox}{\vrule height6pt width6pt depth0pt}

\makeatletter
\renewcommand{\subsection}{\@startsection%
{subsection}{2}{0mm}{\baselineskip}{-1em}%
{\normalfont\normalsize\bfseries}}
\makeatother

\newenvironment{demo}{\noindent{\it Proof.}}
{{\unskip\nobreak\hfil\qquad
\demobox\parfillskip=0pt\par}
\medskip}

\input mathdefs

\renewcommand{\P}{\mathbb{P}}
\newcommand{\Img}{\mathrm{Im}}
\newcommand{\Irr}{\mathrm{Irr}}
\newcommand{\GL}{\mathrm{GL}}

\title{Rigidity and Frobenius Structure}
\author{Richard Crew\\ The University of Florida}

\begin{document}
\maketitle

\section*{Introduction}
\label{sec:intro}

The purpose of this note is show that an irreducible rigid
differential equation on an open subset of $\P^1$ with regular
singularities and rational exponents has, with reasonable local
assumptions relative a prime $p$, a Frobenius structure relative to
some power of $p$.

Katz showed that any irreducible rigid local system on an open subset
$\P^1$ can be built up by repeated tensor product and convolution
operations of a suitable sort from local systems of rank one
\cite{katz:1996}. One therefore expects that if the corresponding
regular singular differential equation is defined, say, over $\bQ$ and
has rational exponents, it should have a Frobenius structure for
almost all $p$. What we show in this paper, in effect, is that if the
differential equation has rational exponents, defines an
overconvergent isocrystal for some prime $p$ and satisfies a few other
reasonable local conditions, it will have a Frobenius structure for
that particular $p$. It is well known that the existence of a
Frobenius structure implies that the equation comes from a convergent
isocrystals, and overconvergence then follows from the other
assumptions.  We remark that when the equation is irreducible, this
Frobenius structure is unique up to a scalar multiple, as was shown by
Dwork \cite{dwork:1989}.

Katz used the theory of algebraic $D$-modules in \cite{katz:1996};
Berthelot's theory of arithmetic $D$-modules is not \textit{a priori}
applicable here since it relies heavily on the existence of a
Frobenius structure (it is not known how to define ``holonomic''
without one). On the other hand, once an overconvergent isocrystal is
known to have a Frobenius structure, its direct image by
specialization is to be a holonomic $\Ddag$-module to which the
methods of \cite{katz:1996} could be applied. The present approach is
elementary in that it uses only the cohomological criterion for
rigidity, together with a $p$-adic analogue (theorem 1 below) in terms
of rigid cohomology. The main point is that if a regular singular
differental equation on an open subset of $\P^1$ is rigid and
irreducible, and defines an overconvergent isocrystal, then that
isocrystal is $p$-adically rigid (theorem 2). The existence of a
Frobenius structure follows from this, assuming rational exponents and
other suitable conditions (theorem 3).

\textit{Acknowledgements.} I would like to thank Francesco Baldassarri
for some useful discussions, and Shishir Agrawal for correcting some
misprints in an earlier version of this paper. I am grateful to the
referee for pointing out further misprints and making a number of
helpful suggestions.

\section{Classical and $p$-adic Rigidity}
\label{sec:rigidity}

Let $U$ be a nonempty Zariski open subset of $\P^1_\bC$, with
analytification $U^{an}$. We recall that a local system $V$ on
$U^{an}$ is \textit{rigid} if any other local system on $U^{an}$ with
the same local monodromy as $V$ is isomorphic to $V$.  Denote by
$j:U^{an}\to\P^1$ the natural inclusion, and set $S=\P^1\sm U$. Katz
shows that an irreducible $V$ is rigid if and only if
$H^1(\P^1,j_*End(V))=0$, or equivalently if
$\chi(\P^1,j_*End(V))=2$. That this condition is sufficient is
relatively easy, and we will see that it can be extends to the case of
$p$-adic differential equations.

A $p$-adic analogue of the rigidity condition can be formulated for
the category of overconvergent isocrystals on an open subset $\P^1$
over a $p$-adic base. We will assume that the reader is familiar with
this theory, but it will be useful to recall a few basic
constructions.

Fix a complete discrete valuation ring $\V$ of mixed characteristic
$p$, with fraction field $K$ and residue field $k$. Let $\P^1_\V$ be
the projective line considered as a formal $\V$-scheme, $U\subset\P^1$
a nonempty formal affine subscheme with closed fiber $U_k$. The
complement $S=\P^1\sm U_k$ is then a finite set of points. As usual,
$U^{an}\subset\P^1_K$ will denote the corresponding affinoid space; it
is the same as the tube $]U[=]U_k[$ (c.f. \cite{berthelot:1996}).
Recall that in this setting, an overconvergent isocrystal on $U$ can
be identified with a locally free module with an overconvergent
%
%
connection $(M,\nabla)$ over the dagger-algebra
\begin{equation}
  \label{eq:global-dagger-algebra}
  A^\dagger=\limdir_W\Gamma(W,\O_W)
\end{equation}
where $W$ runs over the directed system of strict neighborhoods of
$U^{an}$, i.e. a rigid-analytic open neighborhoods $W$ of $U^{an}$
such that $\{W,\P^1\sm U^{an}\}$ is an admissible cover of $\P^1$.  We
will usually abbreviate $(M,\nabla)$ by $M$.

If $s$ is a point of $S$ and $W$ is a strict neighborhood of
$U^{an}$, the open set $W\cap]s[$ is isomorphic to a rigid-analytic
annulus, and we denote by $\R_s$ the direct limit
\begin{equation}
  \label{eq:Robba}
  \R_s=\limdir_W\Gamma(W\cap]s[,\O_W)
\end{equation}
of the function algebras of these annuli; this is the \textit{Robba
  ring} at $s$. If $x$ is a local parameter of $\P^1_\V$ at $s$
(i.e. reduces to a local parameter of $\P^1_k$ at $s$) then
$\R_s\simeq\R$ where $\R$ is the ``standard'' Robba ring, i.e. the
ring of formal Laurent series in $x$ converging in some annulus
$r<|x|<$. If $s$ is a point of $S$, the natural inclusions
$W\cap]s[\inj W$ induce injective ring homomorphisms
$\Gamma(W,\O_W)\inj\Gamma(W\cap]s[,\O_W)$, whence a continuous ring
homomorphism $A^\dagger\inj\R_s$ for all $s\in S$. If $(M,\nabla)$ is
an overconvergent isocrystal on $U$, we set
\begin{displaymath}
  M_s=\limdir_W\Gamma(W\cap]s[,M)
\end{displaymath}
which, since $M$ is a a coherent $\O_W$-module, is a $\R_s$-module of
finite presentation. The connection on the $\R_s$-module $M_s$ induced
by $\nabla$ will be denoted $\nabla_s$, and finally the pair
$(M_s,\nabla_s)$ will be also denoted by $M_s$; it is an
``overconvergent isocrystal on $\R_s$'' that represents the mondromy of
$M$ about $s$.

We therefore make the following definition. An overconvergent
isocrystal $M$ on $U$ is \textit{$p$-adically rigid} if it has the
following property: if $N$ is another overconvergent isocrystal on $U$
such that $M_s\simeq N_s$ for all $s\in S$, then $M\simeq N$. As in
the classical case we do not make a definition in the case of curves
of higher genus, or varieties of higher dimension (although for curves
of higher genus, the definition of ``weakly rigid'' extends in an
obvious way).

To formulate a cohomological condition for the $p$-adic rigidity of an
overconvergent isocrystal $(M,\nabla)$, we recall that for a complex
local system $V$ on an open $U\sset\P^1_\bC$, $H^1(\P^1_\bC,j_*V)$ is
the same as the parabolic cohomology $H^1_p(U,V)$, i.e. the image of
the forget supports map $H^1_c(U,V)\to H^1(U,V)$. In fact the long
exact sequences arising from the exact triangles
\begin{equation}
  \label{eq:6-term-sequence-for-local-systems}
  \begin{split}
    j_!V\to j_*V&\to \bigoplus_{s\in S}(j_*V)_s\Xto{+1}\\
    j_*V\to Rj_*V&\to \bigoplus_{s\in S}(R^1j_*V)_s[-1]\Xto{+1}
  \end{split}
\end{equation}
reduce to exact sequences
\begin{equation}
  \label{eq:parabolic-cohomology}
  \begin{split}
    0\to H^0(U,V)\to \bigoplus_{s\in S}(j_*V)_s\to
    H^1_c(U,j_*V)\to H^1(\P^1,j_*V)\to 0\\
    0\to H^1(\P^1,j_*V)\to H^1(U,V)\to
    \bigoplus_{s\in S}(R^1j_*V)_s\to H^2(\P^1,j_*V)\to 0
  \end{split}
\end{equation}
and an isomorphism
\begin{equation}
  \label{eq:supplementary-isomorphism}
  H^2_c(U,V)\simeq H^2(\P^1,j_*V).
\end{equation}
The assertion follows from this, given that $H^1_c(U,V)\to
H^1(U,V)$ is induced by the composite $j_!V\to j_*V\to
Rj_*V$. From the definitions and \ref{eq:supplementary-isomorphism} we
get equalities
\begin{equation}
  \label{eq:chi_c-chi_p}
  \begin{split}
    \chi(\P^1,j_*V)&=\dim H^0(U,V)-\dim H^1(\P^1,j_*V)+\dim H^2_c(U,V)\\
    &=\chi_c(U,V)+\sum_{s\in S}\dim V_s.
  \end{split}
\end{equation}

The $p$-adic analogue is straightforward, using rigid cohomology (see
\cite{berthelot:1996} for the general definition, and \cite{crew:1998}
for the case of an affine curve). The first fact we need is the
existence of a six-term exact sequence 
\begin{equation}
  \label{eq:6-term-seq}
  \begin{split}
    0\to H^0(U,M)\to &\bigoplus_{s\in S}H^0_{DR}(M_s)\to H^1_c(U,M)\to\\
    \Xto{\d}H^1(U,M)\to &\bigoplus_{s\in S}H^1_{DR}(M_s)\to H^2_c(U,M)\to0
  \end{split}
\end{equation}
for any overconvergent isocrystal $M$ on $U$
\cite[9.5.2]{crew:1998}. In \ref{eq:6-term-seq} the ``local
cohomology'' $H^i_{DR}(M_s)$ is just the ordinary de Rham cohomology
of $M_s=(M_s,\nabla_s)$. We then define the parabolic cohomology
$H^1_p(U,M)$ by
\begin{equation}
  \label{eq:H^1_p-def}
  H^1_p(U,M)=\Img(\d:H^1_c(U,M)\to H^1(U,M)).
\end{equation}
From this we see that \ref{eq:6-term-seq} is the $p$-adic analogue of
the result of gluing together the exact sequences
\ref{eq:parabolic-cohomology} at the term $H^1(\P^1,j_*V)$.

When $H^1_p(U,M)$ has finite dimension, we can define the
``parabolic'' Euler characteristic of $M$ by analogy with the first
part of \ref{eq:chi_c-chi_p}
\begin{equation}
  \label{eq:p-adic-chi-parabolic}
  \chi_p(M)=\dim H^0(U,M)-\dim H^1_p(U,V)+\dim
  H^2_c(U,M)
\end{equation}
and from \ref{eq:6-term-seq} and \ref{eq:p-adic-chi-parabolic} we get
the equality
\begin{equation}
  \label{eq:chi_c-and_chi_p-alg}
  \chi_p(U,M)=\chi_c(M)+\sum_{s\in S}\dim H^0_{DR}(M_s)
\end{equation}
analogous to second part of \ref{eq:chi_c-chi_p}.

The space $H^1_p(U,M)$ will of course have finite dimension if either
of $H^1(U,M)$ or $H^1_c(U,M)$. The finite-dimensionality of these
latter spaces depends on the behavior of $M$ at the points of $S$. The
next proposition extends slightly the main result of \cite{crew:1998}:

\begin{prop}\label{prop:finite-dimensionality}
  Let $M$ be an overconvergent isocrystal on $U$. If $H^1_{DR}(M_s)$
  has finite dimension for every $s\in S$, the $K$-vector spaces
  $H^1(U,M)$, $H^1_c(U,M)$ and $H^1_p(U,M)$ have finite dimension and
  there are canonical duality isomorphisms
  \begin{equation}
    \label{eq:duality}
    \begin{split}
      H^i(U,M)^\vee&\simeq H^{2-i}_c(U,M^\vee)\\
      H^1_p(U,M)^\vee&\simeq H^1_p(U,M^\vee)
    \end{split}
  \end{equation}
  for $0\le i\le 2$.
\end{prop}
\begin{demo}
  By theorem 9.5 of \cite{crew:1998} it suffices to show that for all
  $s\in U_k$ the $K$-linear map $\nabla_s:M_s\to M_s\tens\Omega^1$ is
  a strict morphism of topological vector spaces. Since $M_s$ and
  $M_s\tens\Omega^1$ are LF-spaces this follows from the next lemma.
\end{demo}

In fact this is standard but I do not know a convenient
reference:

\begin{lemma}
  Suppose $u:V\to W$ is a continuous map of LF-spaces such that
  $\Coker(u)$ has finite dimension. Then $u$ is strict.
\end{lemma}
\begin{demo}
  There is a subspace $H\subset W$ of finite dimension that is an
  algebraic supplement to $u(V)$. Since $H$ is separated, its topology
  is the unique separated topology of a finite-dimensional vector
  space, and $H\oplus V$ is an LF-space. The natural map
  $f:H\oplus V\to W$ is surjective and therefore open by the open
  mapping theorem \cite[Prop. 8.8]{schneider:2002}. Suppose now
  $A\subset V$ is open; then $H\oplus A\subset H\oplus V$ is open and
  consequently $f(H\oplus A)=H+u(A)$ is open. Since
  $u(V)\cap(H+u(A))=u(A)$, $u(A)$ is open in $u(V)$.
\end{demo}

The condition that $\dim(H^1_{DR}(M_s))<\infty$ in proposition
\ref{prop:finite-dimensionality} is a consequence of the ``NL
property'' of Christol and Mebkhout. The definition is rather involved
and we refer the reader to \cite{christol-mebkhout:2002} and the
references therein. The one consequence of this condition we need is
the following: if as before $M$ is an overconvergent isocrystal of
rank $d$ on $U$ and satisfies condition NL at every point of $S$, then
\begin{equation}
  \label{eq:index-formula}
  \chi_c(U,M)=d\chi_c(U)-\sum_{s\in S}\Irr(M_s)
\end{equation}
where $\Irr(M_s)$ is the irregularity of the isocrystal $M_s$, defined
in \cite{christol-mebkhout:2002}. In particular, $\chi_c(U,M)$ only
depends on $U$, the rank of $M$ and the irregularities $\Irr(M_s)$. 

We can now state:

\begin{thm}\label{thm:p-adic-rigidity}
  Suppose $M$ is an irreducible overconvergent isocrystal on
  $U\subset\P^1$ such that $End(M)$ satisfies condition NL at every
  point of $S$. If $\chi_p(End(M))=2$ then $M$ is $p$-adically
  rigid.
\end{thm}
\begin{demo}
  The argument is the same as in \cite{katz:1996}. Suppose that $N$ is
  an overconvergent isocrystal such that $M_s\simeq N_s$ for all $s\in
  S$; in particular $M$ and $N$ have the same rank. Since
  $Hom(M,N)_s\simeq End(M)_s$ for all $s\in S$, $Hom(M,N)$ satisfies
  condition NL at every $s$. Then it follows from $\chi_p(End(M))=2$
  and the index formula \ref{eq:index-formula} that
  $\chi_p(Hom(M,N))=2$, and therefore
  \begin{displaymath}
    \dim H^0(\P^1,Hom(M,N))+\dim H^2_c(\P^1,Hom(M,N))\ge 2.
  \end{displaymath}
  On the other hand $Hom(M,N)$ and $Hom(N,M)$ are dual, so the duality
  \ref{eq:duality} yields
  \begin{displaymath}
    \dim H^0(\P^1,Hom(M,N))+\dim H^0(\P^1,Hom(N,M))\ge 2.
  \end{displaymath}
  and we conclude that one of $Hom(M,N)$, $Hom(N,M)$ is nonzero. Since
  $M$ and $N$ have the same rank and $M$ is irreducible, we conclude
  that $M\simeq N$.
\end{demo}

We note that since $End(M)$ is canonically self-dual, the
irreducibility of $M$ implies that either $\chi_p(M)=2$ or
$\chi_p(M)\le 0$, so that $\chi_p(M)=2$ in this case is equivalent to
$H^1_p(U,End(M))=0$. As in the classical case we can think of
$\dim H^1_p(U,End(M))$ as the number of ``accessory parameters'' of
$M$ (see \cite{katz:1996}, p. 5). I do not know if there is a converse
to theorem \ref{thm:p-adic-rigidity}, as is the case over $\bC$; it
would be of interest to settle this question.

\section{Comparison Theorems}
\label{sec:comparison}

Suppose $M$ is a module with a connection with regular singularites
on, say, an open subset $U$ of $\P^1_\bQ$, and denote by $V$ the
corresponding local system on $U_\bC^{an}$. The aim of this section is
to show, under a few (necessary) assumptions, that if $V$ is rigid,
the $p$-adic completion of $M$ is $p$-adically rigid (one condition,
obviously, is that this $p$-adic completion defines an overconvergent
isocrystal). We need not, however, restrict ourselves to the case
where $M$ is defined over $\bQ$, or over a number field. In fact, the
condition that $V$ be rigid is essentially an algebraic condition on
$M$:

\begin{lemma}\label{prop:rigid-is-algebraic}
  Suppose $M$ is a module with a connection with regular singularities
  on some open subset of $\P^1_K$, where $K$ is a field of
  characteristic zero embeddable into $\bC$. If the local system
  $(M\tens_{K,\iota}\bC)^{an}$ is rigid for one choice of embedding
  $\iota:K\to\bC$, it is rigid for any other choice.
\end{lemma}
\begin{demo}
  By Katz's criterion, it suffices to show that
  $\chi_p((M\tens_{K,\iota}\bC)^{an})=2$ if and only if $\chi_p(M)=2$
  (with the latter defined, say by algebraic $D$-module theory), but
  this is just a special case of the Riemann-Hilbert correspondence.
\end{demo}

If $K$ is any field of characteristic zero and $M$ is a module with
regular connection on $\P^1_K$, we can say that $M$ is \textit{rigid}
if there is an absolutely finitely generated subfield $K_0\subset K$
over which $M$ has a model $M_0$, and an embedding $\iota:K_0\to\bC$
such that $\iota(M_0)^{an}$ is a rigid local system; this is evidently
independent of the choice of model, and, by the lemma, of $\iota$. We
remark that a model over an absolutely finitely generated subfield
always exists.

Now in fact one could give a purely algebraic definition of rigidity,
analogous to the definition for local systems, and with this
definition one could prove that $\chi_p(M)=2$ implies that $M$ is
rigid. The converse, however, would not be available without the above
comparison lemma, since it requires a transcendental argument.

Suppose now $\V$ is a complete discrete valuation ring of mixed
characteristic $p$, with fraction field $K$ and residue field $k$.
Let $S\subset\P^1_\V$ be a closed subscheme that finite, flat and
integral over $\V$ and set $U=\P^1\setminus S$. The $U$ and $S$ that
appeared in the last section are now $\hat U$ (the $p$-adic
completion) and $S_k$. As before, $U^{an}$ is the affinoid space
associated to $\hat U$. Finally we denote by $U_K$, $S_K$ the fibers
of $U$ and $S$ over $K$. Note that $S_K$ can be identified with a
finite subset of the tube $]S_k[$, and in fact every point of $S_K$ is
contained in exactly one disk $]s[$ with $s\in S_k$.

Suppose now that $(M,\nabla)$ (as before, usually referred to as $M$)
is a coherent $\O_U$-module with (integrable) connection. We denote by
$M_K$ the corresponding module with connection on $U_K$. If the formal
horizontal sections of $M_K$ have radius of convergence equal to $1$
at every point of $U_K$, then $M_K$ defines an overconvergent
isocrystal on $\hat U$ which we denote by $M^\dagger$. We are interested
in comparing various properties of $M_K$ and $M^\dagger$, subject to a
number of assumptions. The first is purely geometrical:

\begin{description}
\item[C1] For all $s\in S_k$, the disk $]s[$ contains exactly one
  point of $S_K$, which is a $K$-rational point.
\end{description}

Thus each disk contains at most one singular point of $M_K$. The
remaining conditions refer specifically to $M$. Recall that
$a\in\bZ_p$ is $p$-adic Liouville if for every positive real $r<1$,
$|a-n|<r^{|n|}$ has infinitely many solutions $n\in\bZ$.

\begin{description}
\item[C2] $M$ defines an overconvergent isocrystal $M^\dagger$ on
  $\hat U$.
\item[C3] $M_K$ is regular singular, and the exponents of
  $End(M_K)$ belong to $\bZ_p$ and are not $p$-adic Liouville numbers
  (in particular the exponents of $M_K$ itself do not differ by
  $p$-adic Liouville numbers).
\end{description}

In the next theorem and further on we will need a consequence of
Christol's transfer theorem \cite[thm. 1]{christol:1984}, which can be
stated as follows. First, if $A$ is any $n\times n$ matrix $A$ with
entries in $K$, we denote by $M_A$ the free $\R$-module $\R^n$ with
connection given by
\begin{equation}\label{eq:standard-log-module}
  \nabla(u)=du+Au\tens\frac{dx}{x}
\end{equation}
where $x$ is the parameter of $\R$. Recall finally that for $s\in
S_k$, a local parameter of $\P^1_\V$ at $s$ fixes an identification
$\R_s=\R$. 

\begin{lemma}\label{lemma:transfer}
  Suppose $(M,\nabla)$ satisfies \textbf{\emph{C2-C3}}. If $s\in S_k$,
  $M^\dagger_s$ is isomorphic as an isocrystal on $\R$ to $M_A$ for
  some $n\times n$ matrix $A$ with entries in $\V$.
\end{lemma}

In fact Christol's theorem is a purely local statement and we refer
the reader to \cite[thm. 3.6]{crew:1996} for an explanation of how the
lemma follows from \cite[thm. 1]{christol:1984}.

\begin{thm}\label{thm:comp}
  Suppose $M$ satisfies conditions \textbf{\emph{C1-C3}}. If $M_K$ is
  irreducible as a module with connection, $M^\dagger$ is irreducible
  as an overconvergent isocrystal. If in addition $M_K$ is rigid,
  $M^\dagger$ is $p$-adically rigid.
\end{thm}
\begin{demo}
  The first part follows from theorem 2.5 of \cite{crew:1996}, which
  asserts that the differential galois group of $M_K$ is isomorphic to
  the differential galois group (in the category of overconvergent
  isocrystals) of $M^\dagger$. Thus if $M_K$ corresponds to an
  irreducible representation of its differential galois group, so does
  $M^\dagger$.

  If $M_K$ is rigid, then $\chi(U_K,j_*End(M_K))=2$, where as before
  $j:U_K\to\P^1_K$ is the inclusion (and $M_K$ is now regarded as a
  local system on $U_K$). By theorem \ref{thm:p-adic-rigidity} it
  suffices to show that $\chi_p(End(M^\dagger))=2$. 

  By \textbf{C3}, $End(M^\dagger)$ satisfies condition NL at every
  point of $S_k$, and furthermore $\Irr_s(End(M^\dagger))=0$ for all
  $s\in S_k$. Thus
  \begin{displaymath}
    \chi_c(End(M^\dagger))=d\chi_c(U)=\chi_c(End(M_K))
  \end{displaymath}
  where $d$ is the rank of $End(M_K)$. One can also deduce this
  equality from the comparison theorem of Baldassarri-Chiarellotto
  \cite{baldassarri-chiarellotto:1994}. 

  To show that $\chi_p(End(M^\dagger))=\chi_p(End(M_K))=2$, it thus
  suffices to show that $(j_*M_K)_s$ and $M^\dagger_s$ have the same
  dimension for all $s\in S_k$. Suppose $t$ is a local parameter of
  $\P^1_\V$ such that $t=0$ defines a point of $S_K$ in $\P^1_K$, and
  its reduction in $\P^1_k$. Then $(j_*M_K)_s$ and $M^\dagger_s$ are the
  spaces horizontal sections of the connection in respectively in the
  ring of formal Laurent series $K((t))$, and in the ring of elements
  of $\R_s$ convergent for $0<|t|<1$. Since the exponents of
  $End(M_K)$ are not $p$-adic Liouville, lemma \ref{lemma:transfer}
  implies that $M^\dagger_s$ is isomorphic, as $\R_s$-module with
  connection, to a free $\R_s$-module with connection given by the
  matrix of 1-forms $A\tens_K dt/t$, where $A$ is a constant
  matrix. The verification that these spaces have the same dimension
  is then straightforward (see \cite[Lemma 3.4]{crew:1996} for the
  case of $M^\dagger$).
\end{demo}

\section{Frobenius Structure}
\label{sec:frobenius}

We now apply theorem \ref{thm:comp} to the question of Frobenius
structures. To the assumptions already made we add:
\begin{description}
\item[C4] The exponents of $M_K$ are rational.
\end{description}
It is known that if $M$
satisfies \textbf{C1-C3} and has a Frobenius structure, then
\textbf{C4} holds as well.

If $q=p^f$ is a power of $p$, we denote by $\phi:\hat U\to\hat U$ a
lifting of the $q^{th}$-power Frobenius of $U_k$. If $t$ is a global
parameter on $\P^1_\V$, we could of course take $\phi(t)=t^q$; the
theorem in this section allows more general choices.  We denote by
$\sigma:\V\to\V$ the restriction of $\phi$.

The main theorem of this section follows from the following
local-to-global principle:

\begin{lemma}\label{lemma:local-to-global}
  Suppose $M_K$ is an irreducible module with connection satisfying
  \textbf{\emph{C3}}. If $M_K$ is rigid and
  $\phi^*M^\dagger_s\simeq M^\dagger_s$ for all $s\in S_k$, $M^\dagger$ has a
  $q^{th}$-power Frobenius structure.
\end{lemma}
\begin{demo}
  By theorem \ref{thm:comp} we know that $M^\dagger$ is irreducible and
  $p$-adically rigid, and the assertion follows from the definition.
\end{demo}

We denote by $N$ the least common multiple of the denominators of the
exponents of $M_K$ at all points of $S_k$. 

\begin{thm}\label{thm:frobenius}
  Suppose $M$ satisfies conditions \textbf{\emph{C1-C4}}. If $M_K$ is
  irreducible and rigid, then $M^\dagger$ has a $q^{th}$-power Frobenius
  structure for any $q=p^f$ such that $q\equiv 1\pmod N$.
\end{thm}

We remarked in the introduction that the Frobenius structure is unique
up to scalar multiples.

\begin{demo}
  Fix an $s\in S_k$ and a local parameter $x$ of the Robba ring
  $\R_s$. By lemma \ref{lemma:transfer} there is an isomorphism
  $M_s\simeq M_A$ (depending on the choice of $x$, i.e.\ on the
  identification $\R_s\simeq\R$) for some $A$ with rational,
  $p$-adically integral eigenvalues. In view of lemma
  \ref{lemma:local-to-global} we must show that $\phi^*M_A\simeq M_A$.

  We first remark that $M_A\simeq M_{qA}$ for $q$ as above. This is
  elementary: we can assume $A$ is in Jordan normal form, since the
  eigenvalues are rational. We reduce immediately to the case when $A$
  is a single Jordan block with eigenvalue $\lambda$; then $qA$ is
  similar to a block with eigenvalue $q\lambda$, say $A'$, and it
  suffices to show that $M_A\simeq M_{A'}$. Since $q\equiv1\pmod N$ we
  can write $q\lambda=\lambda+k$ with $k\in\bZ$, and the map
  $\R^n\to\R^n$ given by $u\mapsto x^ku$ is the desired isomorphism.

  To conclude we show that $\phi^*M_A\simeq M_{qA}$. We give two
  arguments, an elementary one that needs restrictions on $\V$ and a
  second, less elementary one with no restrictions.  \medskip

  \textit{First method.}  Let $\pi$ be a uniformizer of $\V$ and let
  $e$ be its absolute ramification index. For this proof we assume
  that $e<p-1$ (note that this excludes $p=2$). Before going on we
  recall that in general, if $\nabla$ and $\nabla'$ are connections on
  $\R^n$ given by $n\times n$ matrices of 1-forms $B$ and $B'$, then
  an isomorphism $(\R^n,\nabla)\simeq(\R^n,\nabla')$ is a matrix
  $C\in\GL_n(\R)$ such that
  \begin{equation}
    \label{eq:base-change}
    dC\cdot C^{-1}=CBC^{-1}-B'.
  \end{equation}
  In particular if $B$, $B'$ are conjugate by a constant matrix, that
  matrix also yields an isomorphism $(\R^n,\nabla)\simeq(\R^n,\nabla')$.

  Now
  \begin{equation}
    \label{eq:base-change2}
    \phi^*\left(A\tens\frac{dx}{x}\right)
    =A\tens\frac{d\phi(x)}{\phi(x)}
    =qA\tens\frac{dx}{x}+A\tens\frac{d h(x)}{h(x)}    
  \end{equation}
  with $h(x)=x^{-q}\phi(x)$. We need a $C$ satisfying
  \ref{eq:base-change}, where $B$ is the right hand side of
  \ref{eq:base-change2} and $B'=qA\tens dx/x$. We will find one that
  commutes with $B$, in which case \ref{eq:base-change} reduces to
  $dC\cdot C^{-1}=A\tens dh/h$. If we denote by $\R^0$ the integral
  Robba ring, i.e. the subring of $\R$ with coefficients in $\V$, then
  $h(x)\equiv1\pmod{\pi\R^0}$. We may then define $\log h(x)$ by the
  usual power series, and $\log h(x)\equiv1\pmod{\pi\R^0}$ as
  well. Since $e<p-1$ the exponential $C(x)=\exp(A\tens\log h(x))$
  converges to an element of $\GL_n(\R^0)$. Since $C(x)$ commutes with
  $A$, the change of basis by $C(x)$ is the desired isomorphism
  $\phi^*M_A\simeq M_{qA}$.  \medskip

  \textit{Second method.}  Let $t$ be a global parameter on
  $\P^1_\V$. If $s\in S_k$ corresponds to $t=a$ we set $x=t-a$, which
  is a local parameter at $s$. We denote by $\phi_x$ the lifting of
  the $q$th power Frobenius to $\R_s$ defined by $\phi_x(x)=x^q$. For
  this particular lifting it is evident from \ref{eq:base-change} that
  $\phi_x^*M_A\simeq M_{qA}$, so we must show that
  $\phi^*M_s\simeq\phi_x^*M_s$. Because of our choice of $x$, $\phi_x$
  extends to a lifting of Frobenius on all of $\P^1_\V$, namely
  $\phi_x(t)=a^\sigma+(t-a)^q$. Then there is an isomorphism
  $\phi^*M^\dagger\simeq\phi_x^*M^\dagger$ if overconvergent
  isocrystals on $U$; this is the standard ``independence of lifting''
  property of overconvergent isocrystals. More specifically it follows
  from \cite[Prop. 7.1.6]{lestum:2007} with $Y=Y'=\P^1_k$, $X=X'=U$,
  $P=P'=\P^1_\V$, $u_1=\phi$ and $u_2=\phi_x$. Restricting the
  isomorphism $\phi^*M^\dagger\simeq\phi_x^*M^\dagger$ to the tube
  $]s[$ yields $\phi^*M_s\simeq\phi_x^*M_s$.
\end{demo}

\bibliographystyle{amsplain}
\bibliography{rigid}

\end{document}

%% file: mathdefs.tex
\newcommand{\Coker}{{\operatorname{Coker}}}

\newcommand{\limdir}{\varinjlim}

\newcommand{\tens}{\otimes}

\newcommand{\bC}{{\mathbb{C}}}

\newcommand{\bQ}{{\mathbb Q}}

\newcommand{\bZ}{{\mathbb Z}}

\newcommand{\V}{{\cal V}}

\newcommand{\R}{{\cal R}}
\renewcommand{\O}{{\cal O}}

\newcommand{\sm}{\setminus}
\newcommand{\Xto}{\xrightarrow}

\newcommand{\inj}{\hookrightarrow}

\newcommand{\sset}{\subseteq}

\renewcommand{\d}{{\partial}}

\newcommand{\Ddag}{{\cal D}^\dagger}



%% file: rigid.bbl
\providecommand{\bysame}{\leavevmode\hbox to3em{\hrulefill}\thinspace}
\providecommand{\MR}{\relax\ifhmode\unskip\space\fi MR }
\providecommand{\MRhref}[2]{%
  \href{http://www.ams.org/mathscinet-getitem?mr=#1}{#2}
}
\providecommand{\href}[2]{#2}
\begin{thebibliography}{10}

\bibitem{baldassarri-chiarellotto:1994}
{Baldassarri, F. and Chiarellotto, B.}, \emph{{Algebraic versus rigid
  cohomology with logarithmic coefficients}}, {Barsotti Memorial Symposium},
  {Perspectives in Math.}, vol.~15, {Academic Press}, 1994, pp.~11--50.

\bibitem{berthelot:1996}
{Berthelot, Pierre}, \emph{{Cohomologie rigide et cohomologie rigide \`a
  support propre}}, {Preprint IRMAR 96-03}, 1996.

\bibitem{christol:1984}
{Christol, Gilles}, \emph{{Un th\'eor\`eme de transfert pour les disques
  singuli\`eres reguli\`eres}}, {Cohomologie $p$-adique}, {Ast\'erisque}, vol.
  119--120, SMF, 1984, pp.~151--168.

\bibitem{christol-mebkhout:2002}
{Christol, Gilles } and {Mebkhout, Zoghman}, \emph{{Equations diff\'erentielles
  $p$-adiques et coefficients $p$-adiques sur les courbes}}, {Cohomologies
  $p$-adiques et applications arithm\'etiques (II)} ({Pierre Berthelot and
  others}, ed.), {Ast\'erisque}, vol. 279, SMF, 2002, pp.~125--184.

\bibitem{crew:1996}
{Crew, Richard}, \emph{{The differential Galois theory of regular singular
  $p$-adic differential equations}}, {Mathematische Annalen} \textbf{305}
  (1996), 45--64.

\bibitem{crew:1998}
\bysame, \emph{{Finiteness theorems for the cohomology of an overconvergent
  isocrystal on a curve}}, Annales Scientifique de l'Ecole Normale Sup\'erieure
  \textbf{31} (1998), no.~6, 717--763.

\bibitem{dwork:1989}
{Dwork, Bernard}, \emph{{On the Uniqueness of Frobenius Operator on
  Differential Equations}}, {Algebraic Number Theory -- in honor of K.
  Iwasawa}, {Advanced Studies in Pure Mathematics}, vol.~17, 1989, pp.~89--96.

\bibitem{katz:1996}
{Katz, Nicholas M.}, \emph{{Rigid Local Systems}}, {Annals of Math. Studies},
  vol. {139}, {Princeton Univ. Press}, 1996.

\bibitem{lestum:2007}
{Le Stum, Bernard}, \emph{{Rigid Cohomology}}, {Cambridge Univ. Press}, 2007.

\bibitem{schneider:2002}
{Schneider, Peter}, \emph{{Nonarchimedean Functional Analysis}}, {Springer},
  2002.

\end{thebibliography}
